\documentclass[11pt,
]
{article}

\usepackage{amssymb,amstext,amsthm,latexsym}
\usepackage
{amsmath}
\usepackage{amsfonts}
\usepackage{bbm}

\usepackage{mathrsfs}
\usepackage{mathbbol}

\usepackage{mparhack}

\input{sol.sty}

\begin{document}

\title{A generalization of Solovay's $\Sigma$-construction 
with application to intermediate models\thanks
{This study was partially supported by RFBR grant 13-01-00006 
and Caltech.}}

\author{Vladimir Kanovei
}
%


\date{\today}
\maketitle


\begin{abstract}
A \ddd\Sg construction of Solovay is extended to 
the case of intermediate sets which are not necessarily subsets 
of the ground model, with a more transparent description of the 
resulting forcing notion than in the classical paper of 
Grigorieff.
As an application, we prove that, for a given name $t$ 
(not necessarily a name of a subset of the ground model), 
the set of all sets of the form $t[G]$ 
(the $G$-interpretation of $t$),
$G$ being generic over the ground model, is Borel.
This result was first established by Zapletal by 
a descriptive set theoretic argument.%
\end{abstract}

\def\contentsname{}


\punk{Introduction}

A famous \ddd\Sg construction by Solovay \cite{sol} shows
that if $\dP\in\mm$ is a forcing notion in a countable 
transitive model $\mm$, $t\in\mm$ is a \ddd\dP name, and 
$X\sq\mm$ is any set (\eg, a real), then there is a set
$\xha Xt\sq\dP$ such that
\ben
\Renu
\itla{s1}
the inequality 
$\xha Xt\ne\pu$ is necessary and sufficient for there to exist 
a set $G\sq\dP$, \ddd\dP generic over $\mm$ and satisfying 
$X=t[G]$;

\itla{s2}
if a set $G\sq\dP$ is \ddd\dP generic over $\mm$ and $t[G]=X$ 
then 
$G\sq\xha Xt$, $\xha Xt\in\mn X$, and $G$ is 
\ddd{\xha Xt}generic over $\mn X$; 

\itla{s3}
therefore, in \ref{s2},
%
$\mm[G]$ is a \ddd{\xha Xt}generic extension of $\mn X$;

\itla{s4}
in addition, in \ref{s2}, if a set $G\sq\xha Xt$ is 
\ddd{\xha Xt}generic over $\mm$ 
then still $t[G]=X$.
\een

One may ask whether the \ddd\Sg construction 
of Solovay can be
generalized to {\bf arbitrary} sets $X$, not
necessarily those satisfying $X\sq\mm$.
Following common practice, we'll rather write $\ms X$
in this case.

This paper is devoted to this question, and
the main goal will be to define such a 
generalization, 
although on the base of a somewhat more 
complicated auxiliary forcing $\yha Xt$ which 
consists of ``superconditions'', \ie, pairs of the form 
$\ang{p,a}$, where still $p\in\dP$ while $a$ is a finite 
map associating elements of $X$ with their names.
The generalization (Theorem~\ref{tmain})
will be more direct w.\,r.\,t.\ 
\ref{s1} and \ref{s3}, and rather partial w.\,r.\,t.\ 
\ref{s2} and \ref{s4}.

\bre
\lam{fut1}
Note in passing by that if the axiom of choice holds in 
$\ms X$ then a set $x\sq\mm$ can be easily defined in such 
a way that $\mn x=\ms X$, effectively reducing the problem 
to the case $X\sq\mm$ already considered by Solovay; 
therefore our results below will make sense only in 
the case when $\ms X$ is a choiceless model. 
\ere
 
We'll approach the question from a slightly different 
technical standpoint 
than in the classical paper of Grigorieff~\cite{gri} 
where a base for such generalizations was made. 
This 
will allow us to obtain more pointed generalizations. 
For instance, \ref{s3} is obtained in \cite{gri} by a trick 
which involves a collapse forcing on the top of $\dP$ 
(see Section \ref{g}), 
so that the resulting forcing notion in 
\cite{gri} 
has a much less transparent nature than $\yha Xt$ of this 
paper.

As an application, we prove in Section~\ref{pb} that, for 
a given \ddd\dP name $t$, the set of all sets $t[G]$, 
$G\sq\dP$ being generic over $\mm$, is Borel 
(in terms of an appropriate coding of hereditarily countable 
sets by reals).
Immediately, this set is only analytic, of course.  
This result was first established by Zapletal~\cite{t3} by 
a totally different and much less straightforward argument.

\punk{Basic assumptions and definitions}
\las{blan}

\bdf
\lam{bla}
During the course of the paper, we suppose that:\vihm  
\bit
\item[$-$]
$\mm$ is a countable transitive model of $\ZFC$,\vim  

\item[$-$]
$\dP\in\mm$ is a forcing, and 
$p\le q$ means that $p$ is a {\em stronger\/} condition,\vim  

\item[$-$]
$t\in \mm$ is a \ddd\dP name of a transitive set
(so $\dP$ forces ``$t$ is transitive''),\vim  

\item[$-$]\msur
$X$ is a finite or countable transitive set
(not necessarily satisfying conditions 
$X\in \mm$ or $X\sq\mm$).\qed
\eit
\eDf
The assumption of transitivity of $X$ 
does not reduce generality since 
any set $X$ is effectively coded by the transitive set 
$\ans X\cup\tc(X)$. 

Let $\fo$ be $\fo_\dP^\mm$, the \ddd\dP forcing relation
over the ground model $\mm$.

We also assume that a reasonable {\ubf ramified} system of names 
for elements of \ddd\dP generic extensions of $\mm$ is fixed, 
$\dox$ is a canonical name for any $x\in\mm$. 
and $\doG$ is a canonical name for $G$, the generic set. 
If $t$ is a name and $G\sq\dP$ is \ddd\dP generic over $\mm$ 
then let $t[G]$ be the \ddd Ginterpretation of $t$, so that 
$$
\mm[G]=\ens{t[G]}{t\in\mm\text{ is a name}}.
$$

\bdf
\lam{lm}
In this case, if $X=t[G]$ then $\ms X$ 
is the least transitive model of $\ZF$ 
(not necessarily of $\ZFC$) 
containing $X$ (and all sets in the transitive closure of $X$) 
and all sets in the ground model $\mm$.
Obviously $\mm\sq\ms X\sq\mn G$.
\edf 

For any \ddd\dP names $s,t$, we let 
$s\prec t$ mean that $s$ occurs in $t$ as a name of a potential 
element of $t[G]$.
Then the set $\pe t=\ens{s}{s\prec t}$ 
(of all ``potential elements'' of $t$)
belongs to $\mm$ and if $G\sq\dP$ is generic over $\mm$ then
$$
t[G]=\ens{s[G]}{s\in\pe t\land \sus p\in G\,(p\fo s\in t)}\,.
$$
If $d\sq\pe t$ then a condition $p\in\dP$
is called \rit{\ddd dcomplete\/} iff \vtm

1) $p\fo s\in t$ for
all $s\in d$, and \vtm

2) $p$ decides all formulas $s\in s'$ and $s=s'$,
where $s,s'\in d$.\vtm

\noi
If $d$ is infinite then  \ddd dcomplete conditions do not 
necessarily exist.

\punk{\boldmath Superconditions and the set $\Sigma^+$}
\las{suco}

The following definitions introduce the main technical 
instrument used in this paper: superconditions.

\bdf
\lam{p*}
$\pxt Xt$ is the set of all pairs $\ang{p,a}$ such that 
$p\in\dP$, $a$ is a finite partial map, $\dom a\sq\pe t$, 
$\ran a\sq X$, 
$p$ is \ddd{(\dom a)}complete, 
and in addition $a(\dox)=x$ for any $x\in\mm$ such that 
$\dox\in\dom a$.

We order $\pxt Xt$ so that $\ang{p,a}\le\ang{p',a'}$
($\ang{p,a}$ is stronger) 
iff $p\le p'$ in $\dP$ ($p$ is stronger in $\dP$) 
and $a$ extends $a'$ as a function.
\edf

In particular if $p\in\dP$ then $\ang{p,\pu}\in\pxt Xt$.

Pairs in $\pxt Xt$ will be called 
\rit{superconditions}.\snos
{The upper index $+$ will typically denote things 
related to superconditions as opposite to just 
conditions in $\dP$.}
Given a supercondition  $\ang{p,a}\in\pxt Xt$, we'll call 
$p$ its \rit{condition}, and $a$ its \rit{assignment} --- 
because $a$ assigns sets to (some) names 
forced by $p$ to be elements of $t[G]$.

Note: generally 
speaking, superconditions are not members of $\mm$.

We can also observe that the forcing $\pxt Xt$ just defined 
belongs to $\ms X$ and is a subset of the 
product forcing $\dP\ti\coll(\pe t,X)\in\ms X$.

\ble
\lam{qstro}
If\/ $\ang{p,a}\in\pxt Xt$, 
$q\in\dP$, $q\le p$, 
then\/ $\ang{q,a}\in\pxt Xt$.\qed
\ele

Now we define, following Solovay~\cite{sol}, 
a set $\yha Xt$ of all superconditions
$\ang{p,a}$ which, informally speaking, force nothing really
incompatible with the assumption that there is 
a set $G\sq\dP$ generic over $\mm$ and such that $X=t[G]$ and 
$a(s)=s[G]$ for all $s\in\dom a$. 

\bdf
\lam{defu}
We define a set $\zha\ga Xt\sq\pxt Xt$ 
by transfinite induction on $\ga\in\Ord$.
The dependence on $\dP$ in the definition 
is suppressed. 
\bit
\item
$\zha0Xt$ consists of all superconditions $\ang{p,a}\in\pxt Xt$ 
such that if $s,s'\in\dom a$ then $p\fo s\in\text{(or =)}\,s'$ 
iff $a(s)\in\text{(resp., =)}\;\,a(s')$.

\item
If $\ga\in\Ord$ then the set 
$\zha{\ga+1}Xt$ consists of all superconditions
$\ang{p,a}\in\zha{\ga}Xt$ such that 

\quad
-- for any set $D\in\mm\yd D\sq\dP$, dense in $\dP$, 

\quad
-- and any name $s\in\pe t$ and any element $x\in X$, 

\noi
there is a stronger supercondition 
$\ang{q,b}\in \zha{\ga}Xt$ satisfying: 

\quad a) $\ang{q,b}\le\ang{p,a}$ and $q\in D$, 

\quad b) $x\in \ran b$, \ and 
either $s\in\dom b$ or $q\fo s\nin t$.

\item
Finally if $\la$ is a limit ordinal then 
$\zha{\la}Xt=\bigcap_{\ga<\la}\zha{\ga}Xt$.
\eit
The sequence of sets $\zha{\ga}Xt$ is decreasing,
so that there is an ordinal $\la=\la(X,t)$ such that
$\zha{\la+1}Xt=\zha{\la}Xt$; let $\yha Xt=\zha{\la}Xt$.
%
\edf


\ble
\lam{sgu}
If\/ $\ang{p,a}\in\yha Xt$, a set\/ $D\in\mm\yd D\sq\dP$ is
dense in\/ $\dP$, and\/ $s\in\pe t$, $x\in X,$ 
then there is a pair\/ $\ang{q,b}\in \yha Xt$ satisfying$:$ 
$\ang{q,b}\le\ang{p,a}$, $q\in D$, $x\in \ran b$,
and either\/ $s\in\dom b$ or\/ $q\fo s\nin t$.
\ele
\bpf
This holds by definition, 
as $\yha Xt=\zha{\la}Xt=\zha{\la+1}Xt$.
\epf

The next lemma shows that the set $\yha Xt$ is closed under 
weakening. 


\ble
\lam{qweak}
Assume that\/ $\ang{p,a}\in\yha Xt$. 
Then
\ben
\renu
\itla{qweak1} 
if\/ $\ang{q,b}\in\zha0Xt$ and\/ $\ang{p,a}\le \ang{q,b}$ 
then\/ $\ang{q,b}\in\yha Xt\;;$  

\itla{qweak2} 
if\/ $q\in\dP$, $q\ge p$, 
but still\/ 
$\ang{q,a}\in\pxt Xt$, then\/ $\ang{q,a}\in\yha Xt\;.$ 
\een
\ele
\bpf
\ref{qweak1} 
Prove that $\ang{q,b}\in\zha\ga Xt$ by induction on $\ga$. 
The case $\ga=0$ and the limit case are rather obvious. 
Consider the step $\ga\to\ga+1$. 
By the inductive hypothesis, $\ang{q,b}\in\zha\ga Xt$. 
Let $D\in\mm$, $D\sq\dP$ be dense in $\dP$, 
$s\in\pe t$, and $x\in X$.  
As $\ang{p,a}\in\zha{\ga+1}Xt$, by definition 
there is a stronger supercondition 
$\ang{r,c}\in \zha{\ga}Xt$ satisfying: 
$\ang{r,c}\le\ang{p,a}$, $r\in D$, 
$x\in \ran c$, and either $s\in\dom b$ or $r\fo s\nin t$. 
But then $\ang{r,c}\le\ang{q,b}$ as well, and hence 
$\ang{r,c}$ witnesses $\ang{q,b}\in\zha{\ga+1}Xt$.

\ref{qweak2} 
It follows from $\ang{q,a}\in\pxt Xt$ that 
$\ang{q,a}$ belongs to $\zha0Xt$ together with $\ang{p,a}$. 
It remains to refer to \ref{qweak1}.
\epf

We do not claim that if $\ang{p,a}\in\yha Xt$ and 
$q\in\dP\yt q\le p$ is a stronger condition then, 
similarly to Lemma~\ref{qstro},
$\ang{q,a}\in\yha Xt$.
In fact this hardly can be expected, as 
$q$ may strengthen $p$ in wrong way, that is, 
by forcing about $t$ something 
incompatible with the assignment $a$.
Nevertheless, appropriate extensions of superconditions are 
always possible by Lemma \ref{sgu}.

\punk{The main result}
\las{main}

To formulate the main result, we need one more
definition.

\bdf
\lam{gcomp}
If $G\sq\dP$ is a \ddd\dP generic set over\/ $\mm$ then
let $a[G]$ be the function defined on the set\/ 
$\dom{a[G]}=\pep tG=\ens{s\in\pe t}{s[G]\in t[G]}$ 
so that $a[G](s)=s[G]$ for all $s\in\pep tG$. 

If $\Ga\sq\pxt Xt$ then let\vim 
$$
\bay{rcll}
\prj\Ga&=&\ens{p\in\dP}{\sus a\,(\ang{p,a}\in \Ga)} 
& \text{(the \rit{projection} of $\Ga$ onto $\dP$);}\\[1ex]
A[\Ga] &=& \ens{a}{\sus p\,(\ang{p,a}\in \Ga)} 
& \text{(all assignments which occur in $\Ga$);}\\[1ex]
a[\Ga] &=& \bigcup A[\Ga] 
& \text{(the union of assignments in $\Ga$).}\qed
\eay
$$
\eDf

\ble
\lam{aG}
If\/ $G\sq\dP$ is\/ \ddd\dP generic over\/ $\mm$  
then\/ $\ran{a[G]}=t[G]$.\qed
\ele

Now the main theorem follows;
we prove it in the next two sections.

\bte
[compare with claims \ref{s1}, \ref{s2}, \ref{s3}, \ref{s4} in
the introduction]
\lam{tmain}
In the assumptions of Definition~\ref{bla}
the following holds$:$
\ben
\renu
\itla{ss1}
the inequality 
$\yha Xt\ne\pu$ is necessary and sufficient for there to exist 
a set $G\sq\dP$, \ddd\dP generic over $\mm$ and satisfying 
$X=t[G]\;;$

\itla{ss2}
if a set $G\sq\dP$ is \ddd\dP generic over $\mm$ and $t[G]=X$ 
then the set\/
$$ 
G^+=\ens{\ang{p,a}\in\yha Xt}{p\in G\land a\su a[G]}
$$ 
is\/ \ddd{\yha Xt}generic over\/ $\ms X$, and\/ 
$G=\prj{G^+}\;;$ 

\itla{ss3}
hence, in \ref{ss2},
$\mn G$ is a\/ 
\ddd{\yha Xt}generic extension of\/ $\ms X\;;$ 

\itla{ss4}
%
if a set\/ $\Ga\sq\yha Xt$ is\/ \ddd{\yha Xt}generic over\/
$\ms X$ then the set\/ $H=\prj\Ga\sq\dP$ 
is\/ \ddd\dP generic over\/ $\mm$, $t[H]=X$, and\/ $a[\Ga]=a[H]$.
\een
\ete

\punk{The bounding lemma}
\las{bouL}

Here we prove claim \ref{ss1} of Theorem~\ref{tmain}.
We'll show, in particular, that if indeed $X=t[G]$ for a 
generic set $G\sq\dP$ then the essential length of the 
construction of Definition~\ref{defu} is an ordinal in $\mm$ 
(Lemma~\ref{bou}).

We continue to argue in the assumptions of
Definition~\ref{bla}.

\vyk{
\ble
\lam{uv}
If superconditions\/ $\ang{p,u}$ and\/ $\ang{p,v}$ 
(with the same\/ $p$) belong to\/ $\Sg(X)$ then\/ $u,v$ 
are compatible in the sense that\/ $u(s)=v(s)$ for all\/ 
$s\in\dom u\cap\dom v$.
\ele
\bpf
We argue by induction on the rank of $s$ in $\pe t$. 

If the rank is $0$ then $s$ is the name of the empty set. 
We claim that then $u(s)=v(s)=\pu$. 
Indeed suppose that $u(s)=x\ne\pu$ and let $y\in x$; 
then $y\in X$ 
by transitivity. 
By Lemma~\ref{sgu} there is a stronger supercondition 
$\ang{p',u'}\in\Sg(X)$ with $y\in\ran u'$, so that $y=u'(s')$, 
where $s'\in\dom{u'}\sq\pe t$. 
Thus $u'(s')\in u(s)=u'(s)$, so that by definition 
$q\fo s'\in s$, 
which contradicts the assumption that $s$ 
is the name of the empty  set.

Now the step. 
Suppose towards the contrary that $u(s)\ne v(s)$, 
say $y\in u(s)\bez v(s)$.
By Lemma~\ref{sgu} there is a pair of superconditions 
$\ang{p',u'}\yd\ang{p',v'}\in\Sg(X)$ (with the same $p'$) 
and names $s'\in \dom u'\yd t'\in\dom v'$ such that 
$u'(s')=v'(t')=y$. 
Moreover it follows from the hierarchical organization 
of names  and still Lemma~\ref{sgu} 
\epf

}

The next lemma needs some work.

\ble
\lam{sgG}
Assume that\/ $G\sq\dP$ is a\/ \ddd\dP generic set
over\/ $\mm$, and\/ $t[G]=X$.
If\/ $\ang{p,a}\in\pxt Xt$, $a\sq a[G]$, 
and\/ $p\in G$,  then\/ $\ang{p,a}\in\yha Xt$.

In particular, if\/ $p\in G$ then\/ $\ang{p,\pu}\in\yha Xt$.
\ele
\bpf
Prove $\ang{p,a}\in\zha\ga Xt$ by induction on $\ga$.

Assume that $\ga=0$.
By the \ddd{(\dom a)}completeness, 
if $s,s'\in \dom a$ then $p$ decides $s\in s'$.
If $p\fo s\in s'$ then $s[G]\in s'[G]$,
therefore, as $a\sq a[G]$, we have 
$a(s)\in a(s')$. 
Similarly, if $p\fo s\nin s'$ then $a(s)\nin a(s')$.

The step $\ga\to\ga+1$.
Suppose, towards the contrary, that 
$\ang{p,a}\nin\zha{\ga+1}Xt$ but 
$p\in\zha\ga Xt$ by the inductive hypothesis.
By definition, there exist: a set $D\in\mm\yd D\sq\dP$, 
dense in $\dP$, and elements 
$s\in\pe t$, $x\in X$, such that no supercondition
$\ang{q,b}\in \zha\ga Xt$ satisfies all of 
\bce
$\ang{q,b}\le\ang{p,a}$, \
$q\in D$, \ $x\in \ran b$, \ 
and either $s\in\dom b$ or $q\fo s\nin t$.
\ece
By the genericity, there is a condition $q\in G\cap D$, $q\le p$.
As $t[G]=X$, there is a finite assignment 
$b:(\dom b\sq \pe t)\to X$ 
such that 
\begin{quote}
$a\sq b$, $x\in\ran b$, 
$r[G]\in t[G]$ and $b(r)=r[G]$ for every name $r\in\dom b$,
and either $s[G]\nin t[G]$ or $s\in\dom b$. 
\end{quote}
There is a stronger condition $q'\in G\cap D$ such that 
if in fact $s[G]\nin t[G]$ then $q'\fo s\nin t$, and even more, 
$q'$ is \ddd{(\dom b)}complete.
Then $\ang{q',b}\in\zha\ga Xt$ by the inductive hypothesis, 
a contradiction.

The limit step is obvious.
\epf

\ble
\lam{gins}
If\/ $\ang{p,a}\in\yha Xt$ then there is a set\/ $G\sq\dP$, 
\ddd\dP generic over\/ $\mm$, 
and such that\/ $p\in G$ and\/ $t[G]=X$. 
\ele
\bpf
Both the model $\mm$ and the set $X$ are countable; 
therefore  
Lemma~\ref{sgu} allows 
to define a 
decreasing sequence of superconditions $\ang{p_n,a_n}\in\yha Xt$,%
\pagebreak[0]%
$$
\ang{p,u}=\ang{p_0,a_0}\ge\ang{p_1,a_1}\ge\ang{p_2,a_2}\ge\dots\,,
$$
such that the sequence $\sis{p_n}{n\in\om}$ intersects every 
set $D\in\mm\yd D\sq\dP$, dense in $\dP$ --- hence it 
extends to a generic set $G=\ens{p\in\dP}{\sus n\,(p_n\le p)}$, 
and in addition, the union $\vpi=\bigcup_na_n:\dom\vpi\to X$ 
of all assignments $a_n$ satisfies: 
\ben
\aenu
\itla{gins1} 
$\ran\vpi=X$, 
$\dom\vpi\sq\pe t$, \ and 

\itla{gins2} 
for any $s\in\pe t$\,: 

\quad either $s\in\dom\vpi$ --- then $s[G]\in t[G]$, 

\quad or $q\fo s\nin t$ for 
some $q\in G$ --- then $s[G]\nin t[G]$.
\een
Due to the transitivity of both sets 
$t[G]=\ens{s[G]}{s\in\dom\vpi}$ and $X=\ran\vpi$, 
to prove that $t[G]=X$, it suffices to check that 
$\vpi(s)\in \vpi(s')$ iff $s[G]\in s'[G]$, 
for all names $s,s'\in\dom\vpi$. 
By the construction of $\vpi$, there is an index $n$ 
such that $s,s'\in\dom a_n$. 
By definition, condition $p_n\in G$ is \ddd{(\dom a_n)}complete, 
so $p_n$ decides $s\in s'$. 

If $p_n\fo s\in s'$ then $s[G]\in s'[G]$, and on the other hand, 
as $\ang{p_n,a_n}\in\zha0Xt$, we have 
$\vpi(s)=a_n(s)\in a_n(s')=\vpi(s')$.

Similarly, if $q\fo s\nin s'$ then $s[G]\nin s'[G]$ and 
$\vpi(s)\nin \vpi(s')$.
\epf

The next lemma shows that the ordinals $\la(X,t)$ 
as in Definition~\ref{defu} are bounded in $\mm$ whenever 
$\yha Xt\ne\pu$.

\ble
[the bounding lemma]
\lam{bou}
There is an ordinal\/ $\laa(t)\in\mm$ such that\/ 
$\la(t[G],t)<\laa(t)$ for every set\/ $G\sq\dP$,\/ \ddd\dP generic 
over\/ $\mm$. 
Therefore if\/ $G\sq\dP$ is\/ \ddd\dP generic over\/ $\mm$ then\/ 
$\yha{t[G]}t\in\mm$.
\ele
\bpf
Assume that a set $G\sq\dP$ is \ddd\dP generic over  $\mm$. 
Then $X=t[G]\in\mm[G]$, and hence $\la(X,t)$ is an ordinal in 
$\mm$, and its value is forced, over $\mm$ by a condition $p\in G$, 
to be equal to a 
certain ordinal $\la_p(t)\in\mm$.
We let $\laa(t)=\tsup_{p\in\dP}\la_p(t)$. 
The second part of the lemma follows from the first claim since 
$\yha{X}t$ is the result of a straightforward 
absolute inductive construction of length $\laa(t)\in\mm$.
\epf

\bcor
[= claim \ref{ss1} of Theorem~\ref{tmain}]
\lam{sf}
Tfae\/$:$
\ben
\renu
\itla{sf1}
there is a set\/ $G\sq\dP$, \ddd\dP generic over\/ $\mm$, 
such that\/ $t[G]=X\;;$


\itla{sf2}
$\yha Xt\ne\pu\;;$

\itla{sf3}
$\zha{\laa(t)}Xt=\zha{\laa(t)+1}Xt\ne\pu$.
\een
\ecor
\bpf
Use Lemmas~\ref{sgG}, \ref{gins}, \ref{bou}.
\epf

\punk{Intermediate extensions: proof of the main theorem}
\las{interm}

\vyk{
Here we prove (Theorem~\ref{inT}) 
that the whole extension $\mn G$ is a generic 
extension of any submodel of the form $\ms X,$ where 
$X\in\mn G$, and also establish 
(Lemmas \ref{gag} and \ref{gga}) some analogs of 
claims \ref{s2}, \ref{s4} in Introduction in the process. 
The result makes sense only assuming that 
$\ms X$ is a model of $\ZF$ and not necessarily of $\ZFC$, 
see Remark~\ref{fut1}.
}

In continuation of the proof of Theorem~\ref{tmain}, we prove
here claims \ref{ss2}, \ref{ss3}, \ref{ss4} of the theorem.
We continue to argue in the assumptions of  Definition~\ref{bla}.

\ble
[= claim \ref{ss4} of Theorem~\ref{tmain}]
\lam{gag}
If a set\/ $\Ga\sq\yha Xt$ is\/ \ddd{\yha Xt}generic
over\/ $\ms X$ then the set\/ $H=\prj\Ga\sq\dP$ 
is\/ \ddd\dP generic over\/ $\mm$, $t[H]=X$, and\/ $a[\Ga]=a[H]$.
\ele
\bpf 
By Lemma~\ref{sgu}, if a set $D\in\mm$, $D\sq\dP$, is dense 
in $\dP$ then the set 
$D^*=\ens{\ang{p,a}\in\yha Xt}{p\in D}$ is dense in $\yha Xt$ 
and belongs to $\ms X$. 
It follows that $H$ is indeed generic. 

Further, if $\ang{p,a}\in\Ga\sq\yha Xt$ then by definition 
$\dom a\sq\pe t$ is a finite set and if $s\in\dom a$ then 
$p\fo s\in t$ --- hence, as $p\in H$, we have $s[H]\in t[H]$, 
that is, $s\in\pep tH$. 
On the other hand, if $s\in\pep tH$ and $x\in X$ 
then by Lemma~\ref{sgu} 
there is a supercondition $\ang{q,b}\in\Ga$ such that 
$s\in\dom b$ and $x\in\ran b$. 
Therefore $a[\Ga]$ maps $\pep tH$ onto $X$. 

Still by definition, if $\ang{p,a}\in\Ga$ 
and $s,s'\in \dom a$, then $p$ decides both formulas 
$s\in s'$ and $s=s'$, and 
$p\fo s\in s'$ iff $a(s)\in a(s')$, and the same for $=$. 
Therefore, if $s,s'\in\pep tH$ then we have $s[H]=s'[H]$ 
if and only if  $a[\Ga](s)=a[\Ga](s')$. 
We conclude that $a[\Ga]=a[H]$. 

Finally, $t[H]=\ran{a[H]}=\ran{a[\Ga]}=X$. 
\epf

\ble
[= claim \ref{ss2} of Theorem~\ref{tmain}]
\lam{gga}
If a set\/ $G\sq\dP$ is\/ \ddd\dP generic over\/ $\mm$ 
and\/ $X=t[G]$ then the set\/
$$ 
G^+=\ens{\ang{p,a}\in\yha Xt}
{p\in G\land a\su a[G]}
$$ 
is\/ \ddd{\yha Xt}generic over\/ $\ms X$, and\/ 
$G=\prj{G^+}$\,. 
\ele
\bpf[lemma]
Otherwise there is a condition $p_0\in G$ forcing 
the opposite, so that for any set $H\sq\dP$, 
\ddd\dP generic over $\mm$, if $X=t[H]$ and $p_0\in H$ 
then $H^+$ is not \ddd{\yha Xt}generic over\/ $\ms X$. 
By Lemma~\ref{sgG}, $\ang{p_0,\pu}\in\yha Xt$. 

Consider any set $\Ga\sq{\yha Xt}$, 
\ddd{\yha Xt}generic over $\ms X$ and containing
$\ang{p_0,\pu}$. 
Then $H=\prj{(\Ga)}$ is \ddd\dP generic over 
$\ms X$ and $t[H]=X$ by Lemma~\ref{gag}. 
It remains to prove that $\Ga=H^+,$ that is, given a 
supercondition $\ang{p,a}\in\yha Xt$, we have 
$\ang{p,a}\in\Ga$ iff $p\in H$ and $a\su a[H]$.

If $\ang{p,a}\in\Ga$ then by definition $p\in H=\prj\Ga$ and 
$a\sq a[\Ga]$, but $a[\Ga]=a[H]$ by Lemma~\ref{gag}.

To prove the converse, let  $\ang{p,a}\in\yha Xt$, 
$p\in H$, and $a\su a[H]=a[\Ga]$.
We claim that $\ang{p,a}\in\Ga$.
If $s\in\dom a$ then $a\in\dom{a[\Ga]}$, therefore 
by definition there is a condition 
$\ang{p_s,a_s}\in\Ga$ satisfying $a\in \dom{a_s}$. 
It easily follows that there is a supercondition 
$\ang{q,b}\in\Ga$ satisfying $q\le p$ and $\dom a\sq\dom b$. 
Then in fact $a\su b$ because $a,b\su a[\Ga]$. 
Therefore the supercondition $\ang{q,b}\in\Ga$ is stronger
than $\ang{p,a}\in\yha Xt$. 
We conclude that $\ang{p,a}$ belongs to $\Ga$, too.
\epf

\vyk{

\ble
\lam{ggp}
If a set\/ $G\sq\Sg$ is\/ \ddd{\Sg}generic over\/ 
$\mm[X]$ then it is\/ \ddd\dP generic over\/ $\mm$. 
Conversely, if a set\/ $G\sq\dP$ is\/ \ddd\dP generic 
over\/ $\mm$ and\/ $X=t[G]$ then\/ $G$ is\/ 
\ddd{\Sg}generic over\/ $\mm[X]$. 
\ele
\bpf[lemma]
It suffices to show that if a set $D\sq\dP$, $D\in\mm$ 
is dense in $\dP$ then $D'=D\cap\Sg$ is dense in $\Sg$.
Let $q\in\Sg$, so that $\ang{q,\pu}\in\sgp$.
By Lemma~\ref{sgu} there is a supercondition 
$\ang{r,u}\in\sgp$ such that $r\le q$ and $r\in D$. 
Then by definition $r\in D'$, as required.

To prove the converse, note first 
that $G\sq\Sg$ by Lemma~\ref{sgG}. 
Assume towards the contrary that $G$ is not \ddd{\Sg}generic 
over $\mm[X]$.
Then this is forced by a condition $p_0\in G$, so that 
for any set $H\sq\dP$, \ddd\dP generic over $\mm$, 
if $X=t[H]$ and $p_0\in H$ then $H$ 
is not \ddd\Sg generic over $\mm[X]$.

Now consider any set $H\sq\Sg$, 
\ddd{\Sg}generic over $\mm[X]$ and containing $p_0$. 
Then $H$ is \ddd{\dP}generic over $\mm$ by the above, 
a contradiction. 
\epf
}

\ble
[= claim \ref{ss3} of Theorem~\ref{tmain}]
\lam{inT}
If a set\/ $G\sq\dP$ is\/ \ddd\dP generic over\/ $\mm$, 
and\/ $X=t[G]$,  
then\/ $\mn G$ is a\/ 
\ddd{\yha Xt}generic extension of\/ $\ms X$.
\ele
\bpf
Note that  $\mn G=\ms X[G^+]$ in the 
assumptions of Lemma~\ref{gga}.
\epf

\qeDD{Theorem~\ref{tmain}}

\vyk{

\ref{inT2} 
Apply Lemma~\ref{ggp}.\vom

\ref{inT3} 
As above, the contrary assumption implies that there is a 
condition $p_0\in\Sg$ such that $t[H]\ne X$ for any set 
$H\sq\Sg$, \ddd\Sg generic over $\mm[X]$.
Consider any set $\Ga\sq\sgp$, \ddd\sgp generic over $\mm[X]$ 
and containing $\ang{p_0,\pu}$. 
Then $H=\prj\Ga\sq\Sg$ is \ddd\dP generic over $\mm$, 
$p_0\in H$, and 
$t[H]=X$ by Lemma~\ref{gag}, and in addition $H$ is 
\ddd\Sg generic over $\mm$ by Lemma~\ref{ggp}, contrary 
to the choice of $p_0$.
}

\vyk{
\ble
\lam{=x}
If\/ $G\sq\yg X$ is a\/ \ddd\dP generic set
over\/ $\mm$ then\/ $X=t[G]$.
\ele
\bpf
Because of Lemma~\ref{gins}, it suffices to prove that
$t[G]=t[G']$ for any two \ddd\dP generic sets $G,G'\sq\yg X$.
We check that, even more, $s[G]=s[G']$ for each name
$s\in\pe t\cup\ans t$.

We argue by induction on the name rank of $s$ in $\pe t\cup\ans t$. 
If the rank is $0$ then $s$ is the name of the empty set,
so obviously $s[G]=s[G']=\pu$.

Now the step. 
Suppose towards the contrary that $s[G]\ne s[G']$, 
say $y\in s[G]\bez s[G']$. 
By the construction of names, there is a name $s'\in\pe t$ of the 
rank strictly below the rank of $s$, such that $y=s'[G]$. 
Then $y=s'[G']$ too by the inductive hypothesis. 
\epf
}

\punk{An example}
\las{te}

We still argue in the assumptions of  Definition~\ref{bla}. 
Consider the set
$$
\xha Xt=\prj{\yha Xt}=\ens{p\in\dP}{\ang{p,\pu}\in\yha Xt}\,;
$$
thus $\xha Xt\sq\dP$, and if a set $G\sq\dP$ is \ddd\dP generic  
over $\mm$ and $X=t[G]$ then $G\sq \xha Xt$ by Lemma~\ref{sgG}.
%
Is it true that, similarly to the Solovay claim \ref{s2} 
(Introduction), 
the set $G$ is \ddd{\xha Xt}generic  
over $\ms X$?

The following example easily yields a {\em negative\/} answer. 

\bex
\lam{exam}
Let $\dP$ be the finite-support product of the Cohen forcing; 
a typical condition $p$ in $\dP$ is a map, 
$\dom p\sq\om\ti\om$ is a finite set, and $\ran p\sq\om$. 
Any generic set $G\sq \dP$ forces reals 
$x_n[G]$ such that $x_n[G](i)=r$ iff there is $p\in G$ such 
that $p(n,i)=r$.
We let\/ $\dx n$ be the canonical name for the real  
$x_n[G]=\dx n[G]$, 
and let\/  
$t$ be the name of the set\/ $t[G]=\ens{\dx n[G]}{n\in\om}$. 
In other words, $\ms{t[G]}$ is a well-known symmetric generic 
extension in which $\AC$ fails and $t[G]$ is an infinite 
Dedekind-finite set of reals. 

Sets of the form $t[G]$ are non-transitive, hence, to be in 
compliance with Definition~\ref{bla}, we define the 
transitive closure $U(X)=X\cup U$,
where  
$$
U=\om\cup\ens{\ans{m,n}}{m,n\in\om}  
\cup\ens{\ang{m,n}}{m,n\in\om} 
$$
of any $X\sq\om^\om,$ and accordingly let $t'$ 
be the canonical name of the transitive 
set $t'[G]=\ens{\dx n[G]}{n\in\om}\cup U$. 

As sets in $U$ belong to $\mm$,  
it will be not harmful to identify each $u\in U$ with its 
own canonical name $\dot u$. 
Then $\pe{t'}=\ens{\dx n}{n\in\om}\cup U$.\qed
\eex

\ble
[obvious]
\lam{te1}
If\/ $p\in\dP$ and\/ $n,k,r\in\om$ then\/ 
$p\fo \dx n[G](k)=r$ iff\/ 
$\ang{n,k}\in\dom p$ and\/ $p(n,k)=r$.\qed
\ele

If $X\sq\om^\om$ then 
the set $\pxt {X\cup U}{t'}$ of superconditions 
(Definition~\ref{p*}) 
consists of all pairs 
$\ang{p,a}$ such that $p\in\dP$, $a$ is a map, 
$\dom a\sq\ens{\dx n}{n\in\om}\cup U$ is a finite set,   
$\ran a\sq X\cup U,$ $a(u)=u$ for all $u\in U\cap\dom a$, 
$a(\dx n)\in X$ for all $\dx n\in\dom a$,   
and (the completeness of Definition~\ref{p*}!) 
if a name $\dx n$ and a pair $\ang{k,r}$ ($n,k,r\in\om$) belong 
to $\dom a$ then $p$ decides the formula \lap{$\dx n[G](k)=r$}, 
or equivalently, $\ang{n,k}\in \dom p$. 

Note that $a$ is a bijection for any supercondition $\ang{p,a}$ 
since $\dP$ obviously forces any names 
$s\ne s'$ in $\pe{t'}$ to denote different sets. 

By Definition~\ref{gcomp}, if $G\sq\dP$ is a generic 
set over $\mm$ then a map 
$$
a[G]:\ens{\dx n}{n\in\om}\cup U\onto X\cup U
$$ 
is 
defined by 
$a[G](u)=u$ for all $u\in U$ and  
$a[G](\dx n)=\dx n[G]$ for all $n$. 

Recall that 
$\xha{X\cup U}{t'}=\ens{p\in\dP}{\ang{p,\pu}\in\yha{X\cup U}{t'}}$.

\ble
\lam{te2}
In the case considered, 
if a set\/ $G\sq\dP$ is\/ \ddd\dP generic over\/ $\mm$ 
and\/ $X=t[G]$ then\/ 
\ben
\renu
\itla{te21}
$\xha{X\cup U}{t'}=\dP$, and\/ 

\itla{te22}
$\yha{X\cup U}{t'}$ consists of all superconditions\/ 
$\ang{p,a}\in\pxt {X\cup U}{t'}$  such that 
if both a name\/ $\dx n$ and a pair\/ $\ang{k,r}$  
belong to\/ $\dom a$ then\/ $\ang{n,k}\in \dom p$, and\/ 
$p(n,k)=r$ iff\/ $a(\dx n)(k)=r$.
\een
\ele
\bpf
\ref{te21}
Let $p\in\dP$. 
To prove $p\in \xha{X\cup U}{t'}$, 
it suffices, by Lemma~\ref{sgG}, to define a 
\ddd\dP generic set $G'\sq\dP$ 
such that still $t[G']=X$ and $p\in G'$. 

Let $N=\ens{n}{\sus k\,(\ang{n,k}\in\dom p)}$.
The set $t[G]=\ens{\dx m[G]}{m\in\om}$ is topologically dense 
in $\om^\om,$ 
therefore there is a bijection $\pi:N\to\om$ such that if 
$\ang{n,k}\in\dom p$ (hence $n\in N$) 
then $\dx{\pi(n)}[G](k)=p(n,k)$.

Using the permutation invariance of $\dP$, we 
obtain a generic set $G'\sq\dP$ such that 
$\dx{\pi(n)}[G]=\dx n[G']$ for all $n\in N$, 
still $t[G']=t[G]=X$, and even $x_{m}[G]=x_m[G']$ for all 
but finite $m\in \om$. 
Then $p\in G'$, as required.

\ref{te22}
The proof is similar.
\epf

Thus by \ref{te21} the forcing $\xha{X\cup U}{t'}$ 
coincides with the
given forcing $\dP$ in this case.
But the set $G$ cannot be \ddd\dP generic over $\ms X$,  
basically even over any smaller model $\mn{\dx n[G]}$,
as $X=t[G]=\ens{\dx n[G]}{n\in\om}$. 
This answers in the negative the question above in this
section.

\vyk{
Using similar considerations, we can see that the key 
intermediate forcing $\yha{X\cup U}{t'}$ consists, in this 
case, of all superconditions 
$\ang{p,a}\in\pxt {X\cup U}{t'}$  such that 
if both a name $\dx n$ and a pair $\ang{k,r}$  
belong to $\dom a$ then $\ang{n,k}\in \dom p$ and 
$p(n,k)=r$ iff $a(\dx n)(k)=r$.
}%

Using \ref{te22}, we can prove that $\yha{X\cup U}{t'}$ contains 
a coinitial subset in $\ms X,$ order isomorphic to 
$\bcoll({\ens{\dx n}{n\in\om}},{X})$, 
the bijective collapse forcing which consists of all 
finite partial bijections 
${\ens{\dx n}{n\in\om}}\to {X}$. 

\bcor
\lam{cor2}
In the case considered in this section, 
the whole model\/ $\mn G$ is a\/ 
\ddd{\bcoll({\ens{\dx n}{n\in\om}},{X})}generic extension 
of\/ $\ms X$.\qed
\ecor

Most likely this result has been known since early period of 
forcing, although we are unable to nail a suitable reference.

\punk{Grigorieff's argument}
\las{g}  

To compare our approach with the basic technique of
intermediate models introduced 
in \cite{gri}, we present Grigorieff's proof of the following 
more abstract version of Lemma~\ref{inT}.

\bte
\lam{inTa}
In the assumptions of  Definition~\ref{bla}, 
if a set\/ $G\sq\dP$ is\/ \ddd\dP generic over\/ $\mm$ and\/ 
$X=t[G]$, 
then\/ $\mn G$ is a generic extension of\/ $\ms X$.
\ete
\bpf
Let $\al\in\Ord\cap\mm$ be greater than the von-Neumann 
rank of $X$. 
We put $Y=V_\al\cap\ms X$ (then $X\sq Y$)
and let $H\sq\dC=\text{Coll}(\om,Y)$ be 
generic over $\mn G$.\snos
{It seems that we can define $Y=\tc(X)$ without any harm 
for the ensuing arguments.}
Then $\mn G[H]$ is a generic extension of $\mm$ 
by the two-step iterated forcing theorem, and 
easily there is a real $r$ such that $\ms X[H] = \mn r$. 

Applying Solovay's result \ref{s3} (Introduction)
we conclude that the whole model $\mn G[H]$ is a generic 
extension of $\mn r$. 
But $\mn r=\ms X[H]$ is a generic extension of $\ms X$, 
hence $\mn G[H]$ is a generic
extension of $\ms X$ by the two-step iterated forcing theorem. 

Now, $G \sq \ms X$ and $\ms X\sq \ms X[G] = \mn G\sq \mn G[H]$. 
In other words, $\mn G$ is an intermediate model between 
$\ms X$ as the ground model and $\mn G[H]$ as a generic extension 
of $\ms X$ by the choice of $H$. 
To finish the argument, Grigorieff makes use of the following 
result 
(a part of Theorem 2 in \cite[2.14]{gri}, granted to Solovay), 
with quite a nontrivial proof. 

\ble
\lam{gte}
Let\/ $\dP$ be a forcing in\/ $\mm$, and let\/ $G\sq\dP$ 
be generic over\/ $\mm$. 
If\/ $x\in\mn G$ and\/ $x\sq\mm$, then\/ $\mn x$ is a generic 
extension of\/ $\mm$.
\qed
\ele

Now it suffices to apply the lemma for the models 
$\ms X\sq\mn G\sq\mn G[H]$ in the 
role of the models $\mm\sq\mn x\sq\mn G$ in the lemma.
\epf

It would be interesting, of course, to track down in detail all 
forcing transformations in this proof, to see how the resulting 
forcing is related to  
the forcing directly given by Lemma~\ref{inT}. 
The case considered in Section~\ref{te} would be the most 
elementary one.

\punk{The property of being generic-generated is Borel}
\las{pb}

Another consequence of Lemma~\ref{bou} and other results above 
claims that, in the assumptions of Definition~\ref{bla}, the 
set of all sets of the form $t[G]$, $G\sq\dP$ being generic 
over $\mm$, is Borel in terms of an appropriate coding, of 
all (hereditarily countable) sets of this form, by reals. 
This result was first established by Zapletal 
(Lemma 2.4.4 in \cite{t3}) by 
a totally different argument using advanced technique of 
descriptive set theory.

In order to avoid dealing with coding in general setting, we 
present this result only in the simplest nontrivial 
(= not directly covered by Solovay's original result) 
case when $t$ is a name of a set $t[G]$ which is a  
set of reals, by necessity at most countable.

For a real $y\in\bn$, we let 
$R_y=\ens{(y)_n}{n\in\om}\bez \ans{(y)_0}$, where 
$(y)_n\in\bn$  and $(y)_n(k)=y(2^n(2k+1)-1)$ for all 
$n$ and $k$.
Thus $\ens{R_y}{y\in\bn}$ is the set of all at most 
countable sets $R\sq\bn$ (including the empty set).

\bte
\lam{zapt}
In the assumptions of Definition~\ref{bla}, if\/ $\dP$ forces 
that\/ $t[G]$ is a subset of\/ $\bn$ then the set\/ 
$W$ of all reals\/ $y\in\bn$, 
such that\/ $R_y=t[G]$ for a set\/ $G\sq\dP$ generic over\/ 
$\mm$, 
is Borel.
\ete
\bpf
Let $\vartheta$ be the least ordinal not in $\mm$.
By Corollary~\ref{sf}, for a real $y$ to belong to $W$ 
each of the 
two following conditions is necessary and sufficient:
\ben
\Aenu
\itla{zapt1}
there exist an ordinal $\la<\vartheta$ and 
a sequence of sets 
$\zha\ga Xt\yt \ga\le{\la+1}$, where $X=R_y$, satisfying 
Definition~\ref{defu} and such that 
$\zha\la Xt=\zha{\la+1} Xt\ne \pu$;

\itla{zapt2}
for any ordinal $\la<\vartheta$ and any sequence of sets 
$\zha\ga Xt\yt \ga\le{\la+1}$, where $X=R_y$, satisfying 
Definition~\ref{defu}, if 
$\zha\la Xt=\zha{\la+1} Xt$ then 
$\zha\la Xt\ne \pu$.
\een
Condition \ref{zapt1} provides a $\fs11$ definition of 
the set $W$ 
while condition \ref{zapt2} provides a $\fp11$ definition 
of $W,$  
both relative to a real parameter coding the 
\ddd\in structure 
of $\mm$.
\epf

\end{document}